\documentclass{amsart}
\usepackage{latexsym, amsmath, amssymb, amsthm, mathrsfs}
\usepackage[dvipsnames]{xcolor}
\usepackage{graphicx}
\usepackage[colorlinks=true,linkcolor=blue,urlcolor=blue, citecolor=blue]%
  {hyperref}

\usepackage[shortalphabetic]{amsrefs}
\newcommand{\fan}{\Sigma} 
\usepackage{tikz-cd}
\usetikzlibrary{patterns}
\newtheorem{theorem}{Theorem}[section]
\newtheorem{lemma}[theorem]{Lemma}
\newtheorem{proposition}[theorem]{Proposition}
\newtheorem{corollary}[theorem]{Corollary}

\theoremstyle{definition}
\newtheorem{definition}[theorem]{Definition}
\newtheorem{example}[theorem]{Example}

\theoremstyle{remark}
\newtheorem{remark}[theorem]{Remark}
\newtheorem{question}[theorem]{Question}

\numberwithin{equation}{section}

\newcommand{\generators}{\mathcal{G}}
\newcommand{\toricgraph}{P(b)_\generators}
\newcommand{\mfT}{\mathfrak{T}}
\newcommand{\T}{\mathcal{T}}
\newcommand{\N}{\mathcal{N}}
\newcommand{\I}{\mathcal{I}}
\newcommand{\E}{\mathcal{E}}
\newcommand{\CO}{\mathcal{O}}
\newcommand{\OO}{\CO}
\newcommand{\LL}{\OO(D)} 
\newcommand{\torus}{T}
\DeclareMathOperator{\Pic}{Pic}
\DeclareMathOperator{\hilb}{Hilb}
\DeclareMathOperator{\spec}{Spec}

\DeclareMathOperator{\nefcone}{Nef}

\DeclareMathOperator{\Hom}{Hom}
\DeclareMathOperator{\rk}{rk}
\DeclareMathOperator{\Spec}{Spec}

\DeclareMathOperator{\Nef}{\nefcone}

\DeclareMathOperator{\CL}{Cl}

\newcommand{\KK}{\ensuremath{\mathbb{C}}}
\newcommand{\CC}{\KK}
\newcommand{\QQ}{\ensuremath{\mathbb{Q}}}
\newcommand{\ZZ}{\ensuremath{\mathbb{Z}}}
\newcommand{\RR}{\ensuremath{\mathbb{R}}}
\newcommand{\cA}{\ensuremath{\mathcal{A}}}
\newcommand{\A}{\ensuremath{\mathbb{A}}}
\newcommand{\PP}{\ensuremath{\mathbb{P}}}
\renewcommand{\O}{\ensuremath{\mathcal{O}}}

\begin{document}

\title{Algebraic Hyperbolicity for Surfaces in Toric
  Threefolds}

\author[Haase]{Christian Haase}
\address[Christian Haase]{Institut f\"ur Mathematik\\ Freie Universit\"at
  Berlin\\ Arnimallee~3, 14195 Berlin\\ Germany}
\email{haase@math.fu-berlin.de}
\thanks{This project originated during the Fields Institute program on
  Combinatorial Algebraic Geometry. We thank Fields for
  support. Gregory Smith and Zach Teitler were involved in many useful
  discussions concerning this project. We also thank Sandra Di Rocco
  for several conversations. Luca Chiantini and Angelo Felice Lopez
  were very helpful in explaining their methods to us. Izzet Coskun, Eric Riedl, and Sharon Robins provided helpful comments on an earlier version of this manuscript.}
\thanks{Work of the first author was partially supported by the grant
  HA 4383/8-1 of the German Research Foundation DFG}

\author[Ilten]{Nathan Ilten}
\address[Nathan Ilten]{Department of Mathematics\\ Simon Fraser University\\ 8888
  University Drive\\ Burnaby BC V5A~1S6\\ Canada}
\email{nilten@sfu.ca}
\thanks{Work of the second author was partially supported by the grant
  346300 for IMPAN from the Simons Foundation and the matching
  2015-2019 Polish MNiSW fund.}


\begin{abstract}
Adapting focal loci techniques used by Chiantini and Lopez, we provide
lower bounds on the genera of curves contained in very general
surfaces in Gorenstein toric threefolds.
We illustrate the utility of these bounds by obtaining results on
algebraic hyperbolicity of very general surfaces in toric threefolds.
\end{abstract}

\maketitle

\section{Introduction}
\subsection{Background}
Let $Y$ be a smooth variety over $\KK$. The variety $Y$ is said to be \emph{algebraically hyperbolic} if there is an ample divisor $H$ on $Y$ and some $\epsilon>0$ such that for every integral curve $C\subset Y$,
\[
2g(C)-2\geq \epsilon (C.H).
\]
Here, $g(C)$ is the geometric genus of $C$.
This has been conjectured by Demailly \cite{demailly} to be equivalent to Brody hyperbolicity, see e.g.~\cite{brody}.

The algebraic hyperbolicity of very general surfaces in $\PP^3$ is now completely understood. 
Xu \cite{xu} improved on results of Ein \cite{ein} to show that for very general $S$ of degree $d$ at least $5$,
\begin{equation}\label{eqn:bound}
	g(C)\geq \frac{(d-5)C.H}{2}+2,
\end{equation}
where $H$ is the hyperplane class on $S$, implying that very general surfaces of degree at least $6$ are algebraically hyperbolic. This was recently further improved by Coskun and Riedl \cite{coskun} to the bound
\begin{equation*}
	g(C)\geq \frac{(d(d-5)+1)C.H}{2d}+1,
\end{equation*}
showing that also a very general quintic surface is algebraically hyperbolic. Surfaces of degree at most four contain rational curves, and hence cannot be algebraically hyperbolic.

Very general surfaces in $\PP^3$ are nonetheless very special; for example, very general surfaces of degree at least four always have Picard number one by the Noether-Lefschetz theorem. In this article, we will expand the study of algebraic hyperbolicity to a much larger class of surfaces, namely, very general surfaces in Gorenstein toric threefolds. 

\subsection{Approach and Results}
Xu's bound \eqref{eqn:bound} was originally obtained using a delicate analysis of the behaviour of curve singularities under deformation. This bound 
was subsequently reproven by Chiantini and Lopez using the techniques of \emph{focal loci} \cite{chiantini-lopez}.
It is these focal loci techniques that we use to generalize the bound \eqref{eqn:bound} to very general surfaces in a Gorenstein toric threefold $X$.

Our strongest result (Theorem \ref{thm:z2}) is a bit technical to formulate, but to give a preview we state the following weaker result here:
\begin{theorem}\label{thm:1}
Let $X$ be a Gorenstein toric threefold with torus $\torus$ and let
$D$ be a very ample divisor giving a projectively normal embedding.
For $m\geq 2$, let $S$ be a very general surface in $|mD|$ and
$C\subset S$ an integral curve that is not contained in the toric
boundary $X \setminus \torus$ of $X$. Then the geometric genus $g$ of
$C$ satisfies
\[
	g\geq \frac{C.((m-1)D+K_X)}{2}+1.
\]
\end{theorem}
Although the above theorem only covers curves not contained in the toric boundary, the surface $S$ contains only finitely many curves contained in the toric boundary. Their genera may be determined using combinatorial methods, see Lemma \ref{lemma:boundary}.
In the case that $X=\PP^3$, we may take $D$ to be the hyperplane section $H$, and $K_X=-4H$. The bound resulting from our Theorem \ref{thm:1} becomes the same as that of \eqref{eqn:bound} with the constant term decreased by one.

Applying our stronger bound (Theorem \ref{thm:z2}) we are able to obtain results on algebraic hyperbolicity.
In general, we can show the following.
\begin{theorem}\label{thm:2}
Let $X$ be a Gorenstein projective toric threefold  with nef cone denoted by $\Nef(X)$. There exists an ample divisor class $H_0$ such that for all divisors $D$ whose class lies in $H_0+\Nef(X)$, a very general surface $S\in |D|$ is algebraically hyperbolic. 
\end{theorem}

For specific toric threefolds, we have stronger results:
\begin{theorem}[See Example \ref{ex:p2p1}]\label{thm:p2p1}
Let $X=\PP^2\times \PP^1$. A very general section $S$ of $\CO(a,b)$ is algebraically hyperbolic if $a\geq5$, $b\geq 3$.
\end{theorem}

\begin{theorem}[See Example \ref{ex:p1p1p1}]\label{thm:p1p1p1}
  Let $X=\PP^1\times\PP^1\times\PP^1$. 
  A very general section $S$ of $\CO(a,b,c)$ is
  algebraically hyperbolic if $a \geq b > c = 3$ or if $a
  \geq b \geq c \geq 4$. 
\end{theorem}

\begin{theorem}[See Example \ref{ex:bl}]\label{thm:blp3}
Let $X$ be the blowup of $\PP^3$ at a point, $H$ the pullback of the hyperplane class, and $E$ the exceptional divisor. The nef cone of $X$ is generated by $H$, $H-E$, and a very general section $S$ of $\CO(aH+b(H-E))$ is algebraically hyperbolic if $a\geq 3$, $b\geq 4$ or $a=2$ and $b\geq 7$.
\end{theorem}

We can even apply our methods in non-Gorenstein cases:
\begin{theorem}[See Example \ref{ex:wps}]\label{thm:wps}
Let $X$ be the weighted projective space $\PP(1,1,1,n)$ and $H$ the ample generator of $\Pic(X)$. A very general section $S$ of $mH$ is algebraically hyperbolic if $n\geq 3,m\geq 4$ or $n=2,m\geq 5$.
\end{theorem}

The bounds we obtain for algebraic hyperbolicity in these theorems are close to sharp, leaving only a few cases unresolved. 
See Section \ref{sec:ex} and Question \ref{q:hyper} for details.

\subsection{Other Related Work}
The arguments of \cites{xu,coskun} make strong use of the Noether-Lefschetz theorem: any curve on a surface $S\subset \PP^3$ of degree at least $5$ is a complete intersection of $S$ with some other divisor. On the other hand, \cite{chiantini-lopez} largely avoids using this fact. We also do not use (a toric analogue of) this fact for our main results. However, Bruzo and Grassi have shown that a Noether-Lefschetz type theorem does hold for \emph{some} toric threefolds, see \cites{nl1, nl2}. We can use this in some situations to obtain better lower bounds on the intersection numbers of our curves $C$ with divisors on $X$.

Our Theorem \ref{thm:2} is related to \cite{brotbek}, which shows that for any smooth projective variety $X$, there exists a number $m_0$ such that for any ample divisor $H$ and all $m\geq m_0$, a general hypersurface in $|mH|$ is Brody hyperbolic. In particular, such a general hypersurface is algebraically hyperbolic. While this result applies in considerable more generality than our Theorem \ref{thm:2}, it does not imply the latter.

While we are giving \emph{lower bounds} on the geometric genus of curves on a very general surface $S$ contained in a toric threefold $X$, one may instead ask for \emph{upper bounds} on the genus. In the situation that toric Noether-Lefschetz applies, such curves are complete intersections in the toric threefold $X$ of $S$ with some other divisor $E$. As long as $E$ is nef and big, the genus of a generic curve of this type can be computed using the method of Danilov and Khovanskii from \cite{danilov-khovanski}; this gives an upper bound on the genus of any curve of this type.

\subsection{Organization}
We now describe the organization of the remainder of this paper. In Section \ref{sec:focal}, we recall basics on the theory of focal sets and adapt several claims from \cite{chiantini-lopez} to our setting. The hard work of the paper is done in Section \ref{sec:main}, where we prove our main technical Theorem \ref{thm:z2}. Section \ref{sec:comb} recalls some of the combinatorics associated with toric varieties and uses this to formulate sufficient conditions for applying Theorem \ref{thm:z2}. We put this all together in Section \ref{sec:proofs} to prove our main Theorems \ref{thm:1} and \ref{thm:2}. Finally, in Section \ref{sec:ex} we consider a number of examples and prove Theorems \ref{thm:p2p1}, \ref{thm:p1p1p1}, \ref{thm:blp3}, and \ref{thm:wps}.



\section{Smooth Families and Focal Sets}\label{sec:focal}
In this section, we adapt the techniques of focal loci for our purposes (see e.g. \cite{cc}).
We let $X$ and $B$ be smooth varieties over $\CC$, with $X$ projective. Consider a subvariety $W\subset B\times X$ flat over $B$ with integral fibers, together with a desingularization $V\to W$. After shrinking $B$, we may assume that the composition $\pi:V\to B$ is a smooth morphism, see \cite[Corollary III.10.7]{hartshorne}.
We thus have the following maps:

\[
  \begin{tikzcd}
    V\arrow[dr]\arrow[drr,"\phi"]\arrow{ddr}[swap]{\pi}\\
    &	W \arrow[hook,r]\arrow[d] & B\times X \arrow[dl,"p_B"]\arrow[r,"p_X"]&X\\
    &	B
  \end{tikzcd}
\]
Here, $p_B$ and $p_X$ are the projections onto $B$ and $X$, respectively.

This induces the following diagram of sheaves on $V$ with exact row and column:
\[
  \begin{tikzcd}
    &&0\arrow[d]\\
    &&\pi^*(\T_B)\arrow[d]\arrow[dr,"\lambda"]\\
    0 \arrow[r] & \T_V\arrow [r] & \phi^* (\T_{B\times X})\arrow [r] \arrow[d]& \N_\phi\arrow[r] & 0\\
    &&(p_X\circ\phi)^*(\T_X)\arrow[d]\\
    &&0\\
  \end{tikzcd}
\]

The exact column arises by pulling back the exact sequence of vector bundles
\[
  \begin{tikzcd}
    0 \arrow[r] & p_B^*(\T_B) \arrow[r] &\T_{B\times X} \arrow[r] & p_X^*(\T_X)\arrow[r] & 0.
  \end{tikzcd}
\]
The sheaf $\N_\phi$ is defined as the cokernel to the differential map $\T_V\to \phi^*(\T_{B\times X})$ and is called the  \emph{normal sheaf to $\phi$}. If $V=W$, that is, $\phi$ is a closed embedding, then $\N_\phi$ is just the normal sheaf for this closed embedding (and is in particular locally free since $V$ and $X$ are smooth, see e.g.~\cite[II.8]{hartshorne}).  In particular, we see that any torsion of $\N_\phi$ is supported on the preimage in $V$ of the singular locus of $W$. 

The map $\lambda:(p_B\circ\phi)^*(\T_B)\to \N_\phi$ is called the \emph{characteristic map} for the family $\pi$. Its rank equals 

\begin{equation}\label{eqn:rank}
\rk \lambda=\dim p_X(W)+\dim B-\dim W,
\end{equation}
see \cite[\S 1]{cc}. Note that the arguments in loc.~cit.~apply verbatim when replacing projective space with our variety $X$.

For any point $\eta\in B$, let $V_\eta$ be the fiber over $\eta$ with $\phi_\eta:V_\eta\to X$ the restriction of $\phi$. Restricting the characteristic map $\lambda$ to this fiber gives the map
\[
	\lambda_\eta:T_{B,\eta}\otimes \O_{V_\eta}\to \N_{\phi_\eta}.
\]
Indeed  
\[
	\pi^*(\T_B)_{|V_\eta}\cong \pi_\eta^*(T_{B,\eta})\cong T_{B,\eta}\otimes \O_{V_\eta}
\]
and $\N_{\phi_\eta}\cong (\N_{\phi})_{V_\eta}$ as in \cite[Proposition 1.4]{cc}.

\begin{remark}\label{rem:firstorder}
	If $V=W$, then $H^0(V_\eta,\N_{\phi_\eta})$ parametrizes first order embedded deformations of $V_\eta$ in $X$, see e.g.~\cite[Theorem 2.4]{deftheory}. In this situation, $\lambda_\eta$ is the map sending a tangent vector of $T_{B,\eta}$ to the section of $H^0(V_\eta,\N_{\phi_\eta})$ associated to the first order deformation obtained by restricting $\pi$ to this tangent direction.

	More generally, if we restrict $\lambda_\eta$ to the open subset of $V_\eta$ avoiding the singular locus of $W$, we locally have a similar description of the map. 
\end{remark}

\begin{example}\label{ex1}
We illustrate the above with a non-projective example.
	Consider $B=\A^1=\spec \KK[t]$, $X=\A^2=\spec \KK[x,y]$, and 
	\[W=V(y-(x-t)^2)\subset \A^2.\]
Then $W=V$ is already smooth over $B$.	
We consider the point $\eta=0 \in B$ and set $R=\KK[x,y]$, $S=R/\langle y-x^2\rangle$. The map $\lambda_\eta$ has the form
\begin{align*}
	T_{B,\eta}\otimes S\cong	S&\to S\cong \Hom_{R}(\langle y-x^2 \rangle,S)\\
	f&\mapsto 2x\cdot f.
\end{align*}
This map has rank one everywhere except at the point $(0,0)$, where the rank drops to $0$. 
\end{example}

For sufficiently generic $\eta \in B$, the rank of $\lambda_\eta$ agrees with that of $\lambda$. This is important for the following definition:
\begin{definition}[cf. {\cite[Definition 2.2]{chiantini-lopez}}]
	Assume that $\rk \lambda_\eta=\rk \lambda$. The \emph{global focal set} $F_\eta$ of the fiber $V_\eta$ consists of all those points $v\in V_\eta$ at which the rank of $\lambda_\eta$ drops, that is, is smaller than
	$\dim p_X(W)+\dim B-\dim W$.\footnote{Our definition of global focal set differs from that of \cite{chiantini-lopez} in that we consider the locus of those points in $V_\eta$ where the rank of $\lambda_\eta$ drops, instead of those points where the rank is smaller than $\dim X+\dim B-\dim W$. However, in the situations where we actually use the focal set (Lemma \ref{lemma:fixed} and Proposition \ref{prop:genus}), our definitions agree.}
\end{definition}

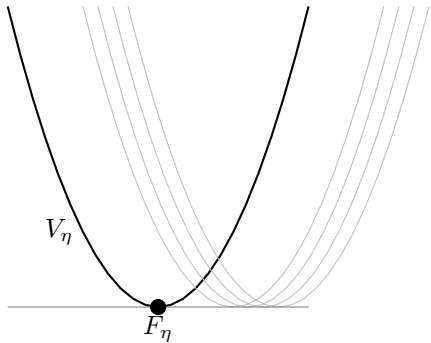
\begin{figure}[htb]
  \centering
\begin{tikzpicture}
\draw[black, thick, domain=-2:2] plot (\x, {\x*\x});
\draw[lightgray,  domain=-2:2] plot (\x+1, {\x*\x});
\draw[lightgray,  domain=-2:2] plot (\x+1.2, {\x*\x});
\draw[lightgray,  domain=-2:2] plot (\x+1.4, {\x*\x});
\draw[lightgray,  domain=-2:2] plot (\x+1.6, {\x*\x});
\draw[gray] (-2,0) -- (2,0);
\draw [black,fill] (0,0) circle [radius=.1];
\node [below] at (0,0) {$F_\eta$};
\node [left] at (-1,1) {$V_\eta$};
\end{tikzpicture}
\caption{A family of parabolas and its focal set (see Example \ref{ex1})}\label{fig:ex1}
\end{figure}

\begin{example}[Example \ref{ex1} continued]
The global focal set in the fiber over $t=0$ of $W=V(y-(x-t)^2)$ consists exactly of the point $(0,0)$. This may be interpreted geometrically as follows. The family $W$ is translating the parabola $V_0=V(y-x^2)$ in the $x$-direction; in particular, the point $(0,0)$ is only contained in the fiber for $t=0$. However, since the line $y=0$ (in the direction of translation) is tangent to $V_0$ at $(0,0)$, this is detected by the global focal set of $V_0$.
See Figure~\ref{fig:ex1}.
\end{example}

For us, the following lemma will play a similar role to	 \cite[Proposition 2.3]{chiantini-lopez}:
\begin{lemma}\label{lemma:fixed}
Suppose that $Y\to B$ is a flat family of hypersurfaces in $X$ such that for every $b\in B$, $W_b\subset Y_b$.
	Let $v\in V_\eta$ be a fixed point of the family $Y$ which avoids the singular loci of $W$, $V_\eta$, and $Y_\eta$. Assume that $\rk \lambda_\eta=\rk \lambda$ and $\dim p_X(W)=\dim X$. Then $v\in F_\eta$ as long as $\rk \lambda>0$.
\end{lemma}
\begin{proof}
Consider some affine chart $\Spec S$ of $X$ containing the image of $v$, and let $I$ be the ideal of $W_\eta$ in this chart. Since $v\in V_\eta$ is mapping to a smooth point of $W$, we may assume that on this chart, $V_\eta=W_\eta$, and the families $V$ and $W$ agree.
Furthermore, since $v$ is a smooth point of both $V_\eta$ and
$Y_\eta$, $V_\eta$ is (locally) a complete intersection in $Y_\eta$,
given say by equations $f_1=\ldots=f_k=0$, where $k=\dim X-1-\dim
V_\eta$. Now, $Y_\eta$ is itself a hypersurface in $X$, cut out by an equation $f_0=0$. 

A tangent vector in $T_{B,\eta}$ determines a first-order embedded deformation $W'$ of $W_\eta$ obtained by pulling back the family $W$; the corresponding element in $\Hom(I,S/I)$ is exactly the image of this tangent vector under $\lambda_\eta$ after restricting to the intersection with the chart $\Spec S$. The family $W'$ (intersected with the chart $\Spec S$) is given by an ideal $\widetilde I$ whose elements are of the form $f+\epsilon f'$ for $f\in I$ and some $f'\in S$. The map $f\mapsto \overline {f'}$ is exactly the element of $\Hom(I,S/I)$ corresponding to this deformation; here $\overline {f'}$ is the image of $f'$ in $S/I$. See e.g.~\cite[Proposition 2.3]{deftheory}. Since the family $W$ is contained in the family $Y$, note that the element $f_0'\in S$ corresponding to $f_0$ may be taken to be the same as the element $f_0'$ defining the pullback of the family $Y$ to the tangent vector of $T_{B,\eta}$.  
Let $J$ be the ideal in $S$ of the point $v$. Restricting $\lambda_\eta$ to $v$ means taking the map $f\mapsto f'$ as an element of $\Hom(I,S/J)$.  
Since $v$ is a fixed point of the family $Y$, we have $f_0'\in J$.
On the other hand, $\Hom(I,S/J)$ is a free module of rank $k+1$, since one can choose the images of $f_1,\ldots,f_k,f_0$ freely. We have seen that under $\lambda_\eta$ restricted to $v$, the image of $f_0$ is always zero. Hence, at $v$, $\lambda_\eta$ can have rank at most $k=\dim X-1-\dim V_\eta$, which is smaller than $\rk \lambda_\eta=\dim X-\dim V_\eta=k+1$.
\end{proof}

We now specialize to the situation where the fibers of $\pi$ are one-dimensional, that is, $V_\eta$ is a smooth curve. Let $g$ denote its genus.
\begin{proposition}\label{prop:genus}
	Assume that $\dim p_X(W)=\dim X$ and the characteristic map $\lambda_\eta$ has the same rank as $\lambda$.
Let $F_\eta^\circ$ be the subset of $F_\eta$ mapping to the smooth locus of $W_\eta$. Then
\[
\deg F_\eta^\circ\leq -K_X.W_\eta+2g-2.
\]
\end{proposition}
\begin{proof}
We adapt the arguments of \cite[Proposition 2.4]{chiantini-lopez}. Using the exact sequence
\[
	\begin{tikzcd}
		0 \arrow[r] & \T_{V_\eta}\arrow [r] & \phi_\eta^* (\T_{X})\arrow [r] & \N_{\phi_\eta}\arrow[r] & 0
\end{tikzcd}
\]
we obtain that 
\[
c_1(\N_{\phi_\eta})=c_1(\phi_\eta^* (\T_{X}))-c_1(\T_{V_\eta})=-K_X.W_\eta+2g-2.
\]

Although the first two terms in the above exact sequence are locally free, $\N_{\phi_\eta}$ might have torsion; let $\mfT$ denote the torsion subsheaf. It is supported on those points of $V_\eta$ mapping to the singular locus of $W_\eta$. 

Consider the composition $\lambda_\eta'$
\[
T_{B,\eta}\otimes \O_{V_\eta}\to \N_{\phi_\eta} \to \N_{\phi_\eta}/\mfT.
\]
The sheaf on the right is now torsion free, hence locally free, since $V_\eta$ is a smooth curve.
 Furthermore, the generic rank of $\lambda_\eta'$ is the same as the rank of $\lambda_\eta$, so the set $F_\eta^\circ$ is contained in the locus where $\lambda_\eta'$ drops rank.

 We now claim that the degree of the locus where $\lambda_\eta'$ drops rank is at most $c_1(\N_{\phi_\eta}/\mfT)$. Indeed, since $\dim p_X(W)=\dim X$, the rank of $\lambda_\eta$ is $\dim X-1$, which is the same as the rank $r$ of $\N_{\phi_\eta}/\mfT$. Hence, $\lambda_\eta'$ is generically surjective. We can thus choose a rank $r$ free subsheaf $\E=\O_{V_\eta}^r$ of $T_{B,\eta}\otimes \O_{V_\eta}$ such that the restriction of $\lambda_\eta'$ to $\E$ still has rank $r$. Now, the locus where this restriction  $(\lambda_\eta')_{|\E}$ drops rank certainly contains the locus where $\lambda_\eta'$ drops rank.

 If $(\lambda_\eta')_{|\E}$ never drops rank, the claim is trivial. Otherwise, the locus where it drops rank has codimension $1=(\rk \E-(r-1))(\rk \N_{\phi_\eta}/\mfT-(r-1))$. This is the ``expected codimension'' of this degeneracy locus, so we can apply the Porteous formula for degeneracy loci (\cite{porteous} or \cite[Corollary 11]{porteousref}). We conclude that the degree of this locus is 
\[
	c_1(\N_{\phi_\eta}/\mfT)/c_1(\E^*)=c_1(\N_{\phi_\eta}/\mfT).
\]

We now obtain that $\deg F_\eta^\circ\leq c_1(\N_{\phi_\eta}/\mfT)$. To conclude, note that $\mfT$ is a direct sum of sheaves supported on a point. Using the standard exact sequence
\[
  \begin{tikzcd}
    0 \arrow[r] & \CO_{V_\eta}(-P)\arrow[r] & \CO_{V_\eta}\arrow [r] & \CO_P\arrow[r] & 0
  \end{tikzcd}
\]
for any point $P\in V_\eta$, we see that a skyscraper sheaf has $c_1=1$. It follows that $c_1(\mfT)\geq 0$, hence $c_1(\N_{\phi_\eta}/\mfT)\leq c_1(\N_{\phi_\eta})=-K_X.W_\eta+2g-2$. 
\end{proof}

\section{Bounding the genus}\label{sec:main}
Let $X$ be a Gorenstein projective toric threefold with torus
$\torus$. We are interested in lower bounds on the genus of curves
contained in general hypersurfaces in $X$. By the \emph{toric
  boundary} of $X$ we mean the complement in $X$ of the open torus
orbit $\torus$.

Let $D$ and $E_1,\dotsc,E_\ell$ be effective, non-trivial torus
invariant divisors on $X$. For each $i$, let $Q_i$ be a basis for
$H^0(X,\CO(E_i))$ consisting of torus equivariant sections. The elements
of each $Q_i$ are uniquely determined up to scaling by a unit of
$\KK$.
\begin{definition}\label{defn:graph}
  The \emph{section graph} for $D,E_1,\ldots,E_\ell$ is the graph $G$ whose vertex set is 
  \[
    V(G)=Q_1\sqcup\cdots\sqcup Q_\ell
  \]
  and where $a\in Q_i$, $b\in Q_{j}$ are connected by an edge if
  and only if there exist $a'\in H^0(X,\CO(D-E_i))$, $b'\in
  H^0(\CO(D-E_{j}))$ such that $aa'=bb'$ in $H^0(\CO(D))$.
\end{definition}

\begin{figure}[htb]
  \centering

\begin{tikzpicture}[y=-1cm]

\draw[black] (13,6.8) circle (0.36cm);
\draw[black] (15,6.8) circle (0.36cm);
\draw[black] (11,6.8) circle (0.36cm);
\draw[black] (11.7,5.6) rectangle (12.3,6.2);
\draw[black] (13.7,5.6) rectangle (14.3,6.2);
\draw[black] (13.7,7.5) rectangle (14.3,8.1);
\draw[black] (11.7,7.5) rectangle (12.3,8.1);
\draw[black] (12.3,5.9) -- (13.7,5.9);
\draw[black] (12.3,7.8) -- (13.7,7.8);
\draw[black] (11.7,7.5) -- (11.25,7.05);
\draw[black] (12.3,7.5) -- (12.75,7.05);
\draw[black] (13.7,7.5) -- (13.25,7.05);
\draw[black] (14.3,7.5) -- (14.75,7.05);
\draw[black] (11.7,6.2) -- (11.3,6.6);
\draw[black] (14.3,6.2) -- (14.7,6.6);
\draw[black] (13.7,6.2) -- (13.3,6.6);
\draw[black] (12.3,6.2) -- (12.7,6.6);
\path (12,7.95) node[text=black,anchor=base] {$1$};
\path (14,7.95) node[text=black,anchor=base] {$x$};
\path (13,6.95) node[text=black,anchor=base] {$x$};
\path (11,6.95) node[text=black,anchor=base] {$1$};
\path (15,6.95) node[text=black,anchor=base] {$x^2$};
\path (12,6) node[text=black,anchor=base] {$y$};
\path (14,6) node[text=black,anchor=base] {$xy$};

\end{tikzpicture}%

  \caption{The section graph for Example \ref{ex:p1p1}. The sections
    of $\CO(E_1)=\CO(1,1)$ are represented by squares, the sections of
    $\CO(E_2)=\CO(1,1)$ by circles.}
  \label{fig:ex:p1p1}
\end{figure}
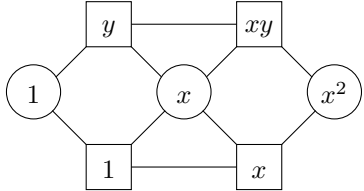

\begin{example}\label{ex:p1p1}
  For illustrative purposes, we consider a section graph for $\PP^1
  \times \PP^1 = (\Spec \KK[x] \cup \{\infty\}) \times (\Spec \KK[y]
  \cup \{\infty\})$ 
  (although it is not a threefold). We consider the configuration
  $D = 2 (\{\infty\} \times \PP^1) + (\PP^1 \times \{\infty\})$
  with global sections $1,x,x^2,y,xy,x^2y$, 
  together with $E_1 = (\{\infty\} \times \PP^1) + (\PP^1 \times \{\infty\})$
  and $E_2 = 2 (\{\infty\} \times \PP^1)$.
  The connected section graph $G$ is pictured in
  Figure~\ref{fig:ex:p1p1}. Observe that $(D;E_1)$ and $(D;E_2)$
  yield the respective induced subgraphs which are not connected.
\end{example}



\begin{definition}\label{defn:cs}
	In the above setting, we say that the configuration of divisors $(D;E_1,\ldots,E_\ell)$ has \emph{connected sections} if 
	\begin{enumerate}
		\item The section graph $G$ is connected;
		\item The union of the images of 
\[H^0(X,\CO(E_i)) \otimes H^0(X,\CO(D-E_i))\]
in $H^0(X,\CO(D))$ span all of $H^0(X,\CO(D))$.
	\end{enumerate}
\end{definition}
\begin{example}
Let $H$ be a torus invariant plane in $\PP^3$. Then for any $d\geq 1$, the configuration $D=dH$, $E_1=(d-1)H$ has connected sections.
\end{example}

In \S\S \ref{sec:markov},\ref{sec:further}, we will give combinatorial criteria for a configuration of divisors to have connected sections. 
\begin{remark}
	Very often, we will only consider a configuration of divisors of the form $(D;E)$, that is, $\ell=1$. However, allowing for $\ell>1$ will give us the flexibility we need to get more refined results, see  Example \ref{ex:cont}.
\end{remark}
\noindent	We now come to our main technical result, from which all other results will follow.
\begin{theorem}\label{thm:z2}
Let $D$ and $E_1,\ldots,E_\ell$ be effective, non-trivial basepoint free torus invariant divisors on a Gorenstein projective toric threefold $X$. Assume that this configuration has connected sections, and that $D$ is big. 
	
Let $S \in |D|$ be a very general surface and $C \subset S$ any integral curve
that is not contained in the toric boundary of $X$.
Then the geometric genus $g$ of $C$ satisfies
\[
  g \geq \min_{i=1,\ldots,\ell} \frac{ C . (E_i + K_{X}) }{ 2 }+1 .
\]
\end{theorem}

\begin{proof}
The first thing that we do is show that we can reduce to the case that $X$ is smooth:
\begin{lemma}\label{lemma:gorenstein}
  Assume that Theorem \ref{thm:z2} is true under the additional
  assumption that $X$ is smooth. Then it also holds for $X$ with
  Gorenstein singularities.
\end{lemma}
\begin{proof}
  Any Gorenstein toric threefold $X$ admits a toric crepant resolution
  $\pi:\widetilde{X}\to X$, see \cite[Proposition 1.1 and
  \S1.2.4]{triangulations}.
  By \emph{crepant}, we mean that $\pi^*(K_X)$ is equivalent to $K_{\widetilde{X}}$.
  We denote the common torus of $X$ and $\widetilde{X}$ by $\torus$.
  If $D,E_1,\ldots,E_\ell$ satisfy the hypotheses of the theorem for
  $X$, then their pullbacks \[\pi^*(D),\pi^*(E_1),\ldots,\pi^*(E_k)\]
  satisfy the hypotheses with respect to the threefold
  $\widetilde{X}$.

  Let $S$ be a very general surface in $|D|$, and $C$ any curve on $S$
  not contained in the toric boundary. Let $\widetilde{S}$ be the
  closure of $S\cap \torus$ in $\widetilde{X}$. Then $\widetilde{S}\in
  |\pi^*(D)|$ is very general. Likewise, let $\widetilde{C}$ be the
  closure of $C\cap\torus$ in $\widetilde{X}$.

  Applying the theorem in the smooth case, we obtain
  \begin{align*}
    g(C)=g(\widetilde{C})
    \geq \min_{i=1,\ldots,\ell} \frac{ \widetilde{C} . \pi^*(E_i + K_{X}) }{ 2 }+1 \\
    = \min_{i=1,\ldots,\ell} \frac{ \pi_*(\widetilde{C}) . (E_i + K_{X}) }{ 2 }+1,
  \end{align*}
  where the second equality follows from the projection formula. But $\pi_*(\widetilde{C})=C$, and the claim follows.
\end{proof}

Using this lemma, we will always assume in the following that $X$ is smooth.
	The idea of the proof, similar to \cite[Theorem 1.3]{cc}, is to construct an appropriate family of curves, show that it has a sufficiently large focal set, and then apply the genus bound from Proposition \ref{prop:genus}.

To begin, fix some number $g$ not satisfying the lower bound of the theorem, and fix numerical invariants $\nu$ determining a Hilbert scheme $\hilb=\hilb_\nu$
of curves in $X$. Let $\hilb'$ be the locus of $\hilb$ parametrizing
integral curves $C$ with geometric genus $g$. This is a locally closed
$\torus$-invariant subscheme. Indeed, the locus of geometrically integral curves is open by \cite[12.2.1(x)]{ega43}. 
Geometric genus is lower semicontinuous by \cite[Proposition 2.4]{harris}, so the subscheme parametrizing curves of fixed geometric genus is locally closed.

The projective space $\PP(H^0(\O(D)))$ parametrizes surfaces in $X$ of class $[D]$.
We now consider the incidence scheme 
\[\I=\left\{ (C,S) \in \hilb' \times \PP(H^0(\O(D)))\ |\ C \subset S
\right\}. \]
This is a closed $\torus$-invariant subscheme of $\hilb'\times
\PP(H^0(\O(D)))$ with projection \[\rho:\I\to \PP(H^0(\O(D))).\]
Let $Z$ be any irreducible component of $\I$, taken with the reduced structure.
Then $Z$ is also invariant under the $\torus$ action.

We will show that as long as $Z$ contains a pair $(C,S)$ with $C$ not contained in the toric boundary of $X$, the image of $Z$ under $\rho$ cannot be dense in $\PP(H^0(\O(D)))$. Then the image of $Z$ under $\rho$ is contained in a proper subvariety of $\PP(H^0(\O(D)))$. The complement of the union of the images of such $Z$ as $\nu$ varies (over a countable set) is a very general subset of $\PP(H^0(\O(D)))$, and by construction, no surface $S$ in this set contains an irreducible curve of genus $g$ not contained in the toric boundary. Hence, the theorem follows.

To show the claim of the previous paragraph, assume that $Z$ contains a pair $(C,S)$ with $C$ not contained in the toric boundary of $X$, and the image of $Z$ under  $\rho$ is dense in $\PP(H^0(\O(D)))$.
Since the image of $\rho$ is constructible (\cite[1.8.4]{ega41})  there is an open subvariety  $U\subset \PP(H^0(\O(D)))$ contained in the image of $\rho$.
Furthermore, $\rho$ has an \'etale section (\cite[17.16.3(ii)]{ega44}), that is, 
we have 
\[
	\begin{tikzcd}
		&Z \arrow[d,"\rho"]\\
		B\arrow[ur,"\sigma"]\arrow[r,"\psi"]&U
	\end{tikzcd}
\]
with $\psi$ \'etale and $\rho\circ \sigma=\psi$. 

Pulling back the universal family of $\hilb'\times\PP(H^0(\O(D))$ along $\sigma$
gives a family $W\subset  B\times X\to B$ of integral genus $g$ curves such that for any $\eta\in B$, the curve $W_\eta$ is contained in the surface $Y_{\eta}$ corresponding to $\psi(\eta)\in \PP(H^0(\O(D))$. 

\begin{lemma}\label{lemma:section}
We may choose the section $\sigma$ such that the image of $W$ in $X$ is dense in $X$.
\end{lemma}
\begin{proof}
Let $T$ be the torus of $X$, which acts linearly on 
$\PP(H^0(\O(D)))$.
Let $H$ be the kernel of this action, and $\overline T=T/H$; this is also a torus, which now acts faithfully on 
$\PP(H^0(\O(D)))$.
Since $D$ is a big divisor, $H$ is a finite subgroup of $T$, so $T\to \overline T$ is \'etale.

After possibly shrinking $U$ (and hence $B$), we can find a subvariety $L\subset U$ such that the rational map $\overline T\times L\dashrightarrow U$ coming from the torus action is birational. 
Indeed, the open torus of $\PP(H^0(\O(D)))$ is a trivial $\overline T$-bundle, so it has such a section $L$.
Furthermore, since a general fiber of the family $W$ is not contained in the toric boundary, we obtain that for a general $\eta\in \psi^{-1}(L)$, $W_\eta$ is also not contained in the toric boundary. 

Take $B'$ to be the locus of $T\times \psi^{-1}(L)$ on which the composition
\[T\times \psi^{-1}(L)\to \overline T\times L\dashrightarrow U\]
is \'etale (and regular). This is non-empty, since  $\overline T\times L\dashrightarrow U$ is birational, and $T\to \overline T$ and $\psi$ are \'etale.
We extend the section $\sigma_{|\psi^{-1}(L)}$ to all of $B'$ via
\[
\sigma'(t,b)=t\cdot \sigma(b)
\]
for $b\in \psi^{-1}(L)$ and $t\in T$. Here we are using the induced $T$-action on $\hilb$ (and thus on $Z$).
We let $W'$ denote the pullback of the universal family of $\hilb'\times\PP(H^0(\O(D)))$ along $\sigma'$. 

By construction, the image of $W'$ in $X$ contains the dense torus, since for general $\eta\in \psi^{-1}(L)$, $W_\eta$  is not contained in the toric boundary.
\end{proof}

We now consider a desingularization $V$ of $W$ (possibly shrinking $B$) so that we are in the situation of \S \ref{sec:focal}. Fix some general point $\eta\in B$.
The surface $Y_\eta$ is the vanishing locus of a section of $\CO(D)$, say $f$. By assumption, the curve $W_\eta$ is not contained in the toric boundary of $X$.
Furthermore, since $X$ is smooth and $D$ is basepoint free, $Y_\eta$
is smooth.

For any surface $S'\subset X$  set
\begin{align*}
B(S')&=\{b\in B\ |\ S'\cap W_\eta\subset Y_b\}.
\end{align*}
This is the preimage under $\psi$ of a linear subspace of $U$; in particular it is smooth.

\begin{lemma}\label{lemma:rank} 
For some $i$ and for a generic $S'\in |E_i|$, the characteristic map for the family $V$ over $B(S')$ has rank $2$ at the point $\eta$.
\end{lemma}
\begin{proof}
	For any section $s\in H^0(X,\CO(E_i))$, let $Y_s^i$ denote the corresponding surface in $X$. Let $G$ be the section graph for the divisors $(D;E_1,\ldots,E_\ell)$ (see Definition \ref{defn:graph}).

Suppose that sections $a\in H^0(X,\CO(E_i))$ and $b\in H^0(X,\CO(E_{j}))$ form an edge in the graph $G$. We will show that the characteristic map at $\eta$ for the base  $B(Y_a^i)\cap B(Y_b^{j})$ has rank at least one. On the other hand, since the union as $i$ varies of the images of $H^0(X,\CO(E_i)) \otimes H^0(X,\CO(D-E_i))$
in $H^0(X,\CO(D))$ span the entire space, it follows that the union of the tangent spaces of $Y_s^i$ at $\eta$ span the tangent space of $B$ at $\eta$. Since $\eta$ was general in $B$ and the dimension of the image of $W$ in $X$ is three, it follows by \eqref{eqn:rank} that the characteristic map at $\eta$ for the base $B$ has rank two (see \S \ref{sec:focal}). We may now apply \cite[Lemma 3.1]{chiantini-lopez} to the characteristic map to conclude that for some subspace $T_\eta B(Y_s^i)$ it has rank two.

It now remains to show that the characteristic map at $\eta$ for the base  $B(Y_a^i)\cap B(Y_b^{j})$ has rank at least one. Let $a'\in H^0(X,\CO(D-E_i))$ and $b'\in H^0(X,\CO(D-E_{j}))$ be such that $aa'=bb'$.
Consider the pencil of surfaces
\[
V(f+t aa')\subset X
\]
(with parameter $t$) corresponding to a line in $\PP(H^0(X,\CO(D)))$ through the point $\psi(\eta)$. Let $A\subset B$ be the preimage of this line under $\psi$. Note that by construction, this curve is contained in $B(Y_a^i)\cap B(Y_b^{j})$ since $aa'=bb'$ and $W_\eta\subset V(f)$. 

Suppose that equations for $W_\eta$ in the torus $T\subset X$ are given by the ideal $I\subset\KK[T]$. Consider a tangent vector of $T_{A,\eta}$ spanning the tangent space; this corresponds to a morphism $\Spec \KK[\epsilon]/\epsilon^2\to A$.
Restricting the family $W$ to the base $\Spec \KK[\epsilon]/\epsilon^2$, the total space intersected with the torus $T$ is cut out by an ideal $\widetilde{I}$ whose elements are of the form 
$h+\epsilon h'$ with $h\in I$ and $h'\in \KK[T]$. 
Away from the singular locus of $W_\eta$ (and after intersecting with the torus), the image of this tangent vector under the characteristic map is the element of $\Hom(I,\KK[T]/I)$ determined by $h\mapsto h'$.
See Remark \ref{rem:firstorder}.

By construction we have $f\in I$. Furthermore, modulo $I$ we must have $f'$ being a non-zero multiple of $aa'$. Then $f'$ is a non-trivial element of $\KK[T]/I$, since otherwise $W_\eta$ would be contained in the zero set of $aa'$, that is, in the toric boundary. It follows that the element of $\Hom(I,\KK[T]/I)$ determined by this tangent vector is non-zero, hence the characteristic map at $\eta$ for the base  $A$ (and hence $B(Y_a^i)\cap B(Y_b^{j})$) has rank at least one. 
\end{proof}

We can now finish the proof of the theorem. Let $S'\in |E_i|$ be as in the lemma above. 
We will now restrict the families $W$ and $Y$ to $B(S')$; by abuse of notation we will still denote them by $W$ and $Y$. The content of the above lemma was that over $B(S')$, the characteristic map for the family $W$ has rank two at $\eta$, which is also its generic rank. Next, we wish to apply Lemma \ref{lemma:fixed}.

We had already noted above that $Y_\eta$ is smooth. Since $E_i$ is basepoint free, we may assume that $S'$ intersects $W_\eta$ transversely and does not contain any point of the singular locus of $W_\eta$.
Furthermore, since we had chosen $\eta$ generically, the family $W$ is smooth at a generic point of $W_\eta$, and this remains true after restricting to the base $B(S')$ \cite[III.10.1]{hartshorne}. Since $B(S')$ is itself smooth at $\eta$, we conclude by loc.~cit.~that $W_\eta$ is not contained in the singular locus of $W$. Thus, we may choose $S'$ so that $S'\cap W_\eta$ does not contain any point of the singular locus of $W$. Finally, we notice that $S'\cap W_\eta$ is contained in every surface $Y_b$  for $b\in B(S')$. By Lemma \ref{lemma:fixed}, we conclude that the global focal set of $V_\eta$ contains $W_\eta\cap S'$.

This implies that the degree of $F_\eta^\circ$ (for the family over $B(S')$) is at least $C.E_i$. On the other hand, since the rank of the characteristic map at $\eta$ has rank $2$, we may apply Proposition \ref{prop:genus} to conclude that
\[
g\geq \frac{C.(E_i+K_X)}{2}+1,
\]
a contradiction. Hence, the image of $Z$ under $\rho$ cannot be dense in $\PP(H^0(\O(D))$, and the theorem is proved.
\end{proof}

\begin{remark}
	In the proof of Theorem \ref{thm:z2}, the fact that $X$ is toric plays a relatively minimal role. We first use the existence of crepant resolutions for toric threefolds to reduce to the smooth case in Lemma \ref{lemma:gorenstein}. We then use the torus action in Lemma \ref{lemma:section} to obtain a section such that the image of the family is dense in $W$. Finally, we use in the proof of Lemma \ref{lemma:rank} that we may work with explicit equations in the coordinate ring of $T$. 

We suspect that the above proof can be adapted to work in other situations where $X$ admits an action with an open orbit by some algebraic group $G$, for example,  when $X$ is an abelian or spherical variety.
\end{remark}

\section{Combinatorial Interpretation and Results}\label{sec:comb}
\subsection{Toric Varieties and Polytopes}\label{sec:toric}
In this section we introduce notation and recall some basic 
facts about toric varieties. For more details 
we refer to~\cite[\S\S2.3,4.1,4.2]{CoxLittleSchenckToricBook}
or~\cite[Section~3.4]{Fulton}. 

Let $N\cong \ZZ^n$ be a lattice with dual
lattice $M = \Hom (N, \ZZ)$ and associated vector spaces $N_{\RR} : =
N\otimes _{\ZZ}\RR$ and $M_{\RR} : = M\otimes _{\ZZ}\RR$.
To a complete rational fan $\fan$ in $N_{\RR}$ we associate 
the toric variety $X = \operatorname{TV}(\fan)$, a normal
equivariant compactification of the algebraic torus $\torus =
\Hom(M,\KK^*)$.
The irreducible components of the toric boundary $X \setminus \torus =
\bigcup_{\rho \in \fan[1]} D_\rho$ are torus invariant prime divisors
$D_\rho$ indexed by the set $\fan[1]$ of one-dimensional cones (rays)
of $\fan$. We can thus identify the group of $\torus$-invariant
Weil divisors with $\ZZ^{\fan[1]}$. It fits into a short exact sequence
(cf.~\cite[Thm.~4.1.3]{CoxLittleSchenckToricBook})
\begin{equation}
  \label{eq:ses-class-gp}
  \begin{tikzcd}
    0 \arrow[r] & M \arrow[r,"\iota"]& \ZZ^{\fan[1]} \arrow[r,"\pi"]& \CL(X) \arrow[r]& 0 \,.
\end{tikzcd}\end{equation}
The $\rho$-coordinate of the map $\iota:M \to \ZZ^{\fan[1]}$ is given by $m
\mapsto \langle v_\rho, m \rangle$ where $v_\rho \in N$ stands for the
primitive generator of the ray $\rho$.
The sequence \eqref{eq:ses-class-gp} implies that every divisor class
contains a $\torus$-invariant divisor $D = \sum_\rho a_\rho
D_\rho$ labeled by $a \in \ZZ^{\fan[1]}$. To the latter, we can
associate the polyhedron
\begin{equation}
  \label{eq:divisor-polytope}
  P = P(D) = \bigl\{ u \in M_{\RR} \mid \langle
  v_\rho, u \rangle + a_\rho \ge 0 \,, \; \rho \in \fan[1] \bigr\}
\end{equation}
which is bounded because we assumed $\fan$ to be complete. Linearly
equivalent invariant divisors yield polytopes which differ by a translation from an element of $M$.
The image of the polytope $P$ under $\iota$ can also be recovered from \eqref{eq:ses-class-gp}:
\begin{equation}\label{eq:polytope}
	\iota(P)=\pi^{-1}([D])_\RR\cap \RR_{\geq 0}^{\fan[1]}- a \,. 
\end{equation}
If $\CL(X)$ is torsion free, then the sequence \eqref{eq:ses-class-gp} splits, so the polytope $P$ is lattice equivalent to $\iota(P)$.

The lattice points in $P$ provide an equivariant basis for the global sections:
\[ H^0(X, \LL) \cong \bigoplus _{u\in P\cap M} \KK \chi^u, \]
where $\chi ^u\colon  \torus \to \KK^*$ is the character corresponding
to $u\in M$.
When $X$ is projective, $D$ is ample if and only if $\fan$ is the
normal fan of $P$.
Further, $D$ is Cartier if there is a continuous function $\varphi_D
\colon N_\RR \to \RR$ given by $m_\sigma \in M$ along cones $\sigma
\in \fan[n]$ so that $a_\rho = \varphi_D(v_\rho)$ for all $\rho$. Then, $D$
is nef if and only if $\varphi_D$ is convex which is equivalent to $P$
being the convex hull of the lattice points $m_\sigma$, $\fan$
refining the normal fan of $P$, and all inequalities in
\eqref{eq:divisor-polytope} being tight. If $D$ and $D'$ are nef, then
\[
P(D+D')=P(D)+P(D'),
\]
where the addition here is the Minkowski sum.
We recall also that for toric varieties, being nef is the same thing
as being basepoint free.

Let now $D$ be a basepoint free Cartier divisor on a toric threefold $X$. The following lemma lets us give a lower bound on the geometric genus of curves contained in the intersection of a general surface in $D$ with the toric boundary:
\begin{lemma}\label{lemma:boundary}
	For $S\in |D|$ a general surface and $C\subset S$ an irreducible curve contained in the toric boundary of $X$, $C=S\cap D_\rho$ for some $\rho\in\fan[1]$ corresponding to a facet $F\prec P(D)$. The geometric genus of $C$ equals the number of interior lattice points of $F$.
\end{lemma}
	\begin{proof}
Since $C$ is irreducible, it must be contained in $S\cap D_\rho$ for some ray $\rho$. The restriction of $D$ to $D_\rho$ is basepoint free, and $D_\rho$ has isolated singularities, so by choosing $S$ general we may assume that $S\cap D_\rho$ is smooth for all $\rho$; it follows that $C=S\cap D_\rho$ and $C$ is smooth.

The statement concerning the geometric genus of $C$ now follows from \cite[Prop.~10.5.8.]{CoxLittleSchenckToricBook}, since $C$ is smooth.
\end{proof}

\subsection{Connected Sections and Markov Bases}\label{sec:markov}
In this section, we will relate the notion of connected sections (Definition \ref{defn:cs}) to the notion of a \emph{Markov basis}.
To that end, we will modify \eqref{eq:ses-class-gp}
to obtain a short exact sequence
\begin{equation}
	\label{eq:modify}
  \begin{tikzcd}
	  0 \arrow[r] & M \arrow[r,"\iota'"]& \ZZ^p \arrow[r,"\pi'"]& \ZZ^q \arrow[r]& 0 \,.
\end{tikzcd}\end{equation}
with the property that for every $T$-invariant Cartier divisor $D=\sum
a_\rho D_\rho$ on $X$, there exists $b\in \ZZ^q$ such that the
polytope $P(D)$ is lattice equivalent, via $\iota'$, to $\RR_{\geq 0}^p\cap (\pi_\RR')^{-1}(b)$.
If $\CL(X)$ is torsion free, 
\eqref{eq:polytope} implies that we may take the sequence \eqref{eq:modify} to just be \eqref{eq:ses-class-gp} after choosing a basis of $\CL(X)$.

For the general case, fix any basis $v_1,\ldots,v_n$ of $N$. We take 
\[
	\ZZ^p=\ZZ^n\times \ZZ^{\fan[1]}
\]
along with the inclusion 
\[
\iota':M\to \ZZ^p\qquad m\mapsto \big((\langle v_i,m\rangle)_i,(\langle v_\rho,m\rangle)_\rho\big).
\]
By construction, this inclusion has a co-section (by projecting to $\ZZ^n$ and using the dual basis). Hence, the cokernel is free, and after choosing a basis we obtain a surjection $\pi':\ZZ^p\to \ZZ^q$.

\begin{lemma}\label{lemma:good}
Let sequence \eqref{eq:modify} be as constructed above. Then for any
$T$-invariant Cartier divisor $D$, there exists $b\in \ZZ^q$
such that $P(D)$ is lattice equivalent to 
\[P(b):=\RR_{\geq 0}^p\cap (\pi_\RR')^{-1}(b).\]
\end{lemma}
\begin{proof}
  The polytope $P(D)$ is defined by inequalities $\langle
  v_\rho,u\rangle \geq -a_\rho$ for $\rho\in\fan[1]$. By setting
  $a'_i := \phi_D(v_i)$, this is the same as imposing the
  inequalities
\[
	\langle v_\rho,u\rangle \geq -a_\rho\qquad\textrm{and}\qquad
	\langle v_i,u\rangle \geq -a'_i
\]
for all $\rho\in\fan[1]$ and $i=1,\ldots,n$. 
Then
\[
\iota'(P(D))=P(b)-(a',a)
\]
where $b=\pi'(a',a)$.
In particular, $P(D)$ is lattice equivalent to $P(b)$ since the
sequence \eqref{eq:modify} is split.
\end{proof}

\begin{remark}
  From a different point of view, we may obtain \eqref{eq:modify} and Lemma \ref{lemma:good} as follows. Let $\phi:\widetilde{X}\to X$ be any toric partial resolution of $X$ so that $\CL(\widetilde{X})$ is torsion free.
  Then the sequence \eqref{eq:modify} may be taken to be the sequence \eqref{eq:ses-class-gp} for $\widetilde{X}$. To obtain the $b\in\ZZ^q\cong \CL(\widetilde X)$ of the lemma for a given Cartier divisor $D$, we take the image in $\CL(\widetilde{X})$ of $\phi^*(D)$. The claim of the lemma follows from \eqref{eq:polytope} and the freeness of $\CL(\widetilde{X})$, along with the fact that $P(D)=P(\phi^*(D))$. 

  A particular instance of such a partial resolution $\widetilde{X}\to
  X$ may be obtained as above: letting $v_1,\ldots,v_n$ be a basis of
  $N$ disjoint from the rays of $\fan$, we consider the stellar
  subdivision of $\Sigma$ along the rays $\rho$ generated by the
  $v_i$. The sequence $\eqref{eq:modify}$ we obtain is exactly the one
  constructed above.
\end{remark}

Having fixed a sequence \eqref{eq:modify} as above, we represent the map 
$\pi':\ZZ^p\to\ZZ^q$ by a matrix $A$.
After choice of bases, the matrix $A$ is determined by $\fan$ alone and does not depend on the divisor $D$.
We are now in the standard situation of \cite[Chapters 4 and 5]{GBCP}
The \emph{toric ideal} $I_A$ associated to the matrix $A$ is the ideal in
$\CC[x_1, \ldots, x_p]$ generated by binomials $x^{v^+}-x^{v^-}$ for
$v^+, v^- \in \ZZ_{\ge 0}^p$ with $Av^+ = Av^-$.

We identify a vector $v \in \ker A \cap \ZZ^p$ with the binomial
\[x^{v^+}-x^{v^-} \in I_A,\]
where $v^+_i = \max(v_i,0)$ and $v^-_i =
-\min(v_i,0)$. Accordingly, we say that a subset $\generators \subset
\ker A \cap \ZZ^p$ is a \emph{Markov basis} if the corresponding binomials
generate $I_A$. For any $b\in \ZZ^q$, we consider a graph
$\toricgraph$ whose vertices are $P(b) \cap \ZZ^p$, and $v,v' \in P(b) \cap \ZZ^p$ are joined by an edge if $v-v' \in
\pm\generators$.

\begin{theorem}[{\cite[Thm~3.1]{DiaconisSturmfels},\cite[Thm~5.3]{GBCP}}
]\label{thm:markov}
  A set $\generators \subset \ker A \cap \ZZ^p$ is a Markov basis for
  the toric ideal $I_A$ if and only if the graph $\toricgraph$ is
  connected for all $b \in \ZZ^q$.
\end{theorem}

We will now apply this theorem to obtain a criterion for connected sections.
Following \cite{convex-normal}, we call a pair $(E,E')$ of nef
divisors \emph{IDP} (it has the integer decomposition property) if 
\[H^0(X,\OO(E)) \otimes H^0(X,\OO(E')) \to H^0(X,\OO(E+E'))\]
is surjective.

\begin{proposition} \label{prop:Markov}
	Let $(E,E')$ be an IDP pair of divisors on $X$, with $b'\in\ZZ^q$ such that $P(E')$ is lattice equivalent to $P(b')$ as in Lemma \ref{lemma:good}.
	Set $D:=E+E'$ and \[\generators := (P(b') \cap \ZZ^p) - (P(b') \cap \ZZ^p).\]
If $\generators$ is a Markov basis for $I_A$, then the configuration $(D;E)$ has connected sections.
\end{proposition}

\begin{proof}
	Since $(E,E')$ is IDP, the second criterion of Definition \ref{defn:cs} is fulfilled. Thus, we only need to show that the section graph for $(D;E)$ is connected.

  The vertices for this graph are in bijection with the lattice points of the polytope $P(E)$; this is lattice equivalent to a polytope $P(b)$ for some $b\in \ZZ^p$.  
  After identifying the vertices of the section graph with the lattice points of $P(b)$, we obtain that the section graph has the same vertex set $P(b) \cap \ZZ^p$ as the toric graph
  $\toricgraph$. The latter is connected by assumption on $\generators$ and Theorem \ref{thm:markov}. We argue that every
  edge of $\toricgraph$ is also an edge in the section graph.

  Indeed, consider $u_1,u_2 \in P(b) \cap \ZZ^p$ corresponding to sections
  \[\chi^{m_1},\chi^{m_2}\in H^0(X,\CO(E)).\] If these vertices are connected by an edge in
  $\toricgraph$, then $u_1-u_2 \in \generators$, that is, there
  are $u'_1,u'_2 \in P(b') \cap \ZZ^p$ corresponding to sections
  $\chi^{m'_1},\chi^{m'_2}$ so that $u_1-u_2 = u'_1-u'_2$.
  On the level of sections, we see that $\chi^{m_1}\chi^{m'_2} =
  \chi^{m_2}\chi^{m'_1} \in H^0(X,\OO(D))$. Hence, $\chi^{m_1}$ and
  $\chi^{m_2}$ are joined by an edge in the section graph.
\end{proof}

\subsection{Further Criteria for Connected Sections}\label{sec:further}
Using the discussion of \S \ref{sec:toric} and \S \ref{sec:markov}, we
will formulate further sufficient criteria for connected sections.
In the following, $X$ is always a projective toric threefold and $H$
an ample divisor.
In dimension three, $(H,H)$ being IDP implies that $(X,H)$ is
projectively normal~\cites{EW91,LTZ93,BGT97}. In particular, for any
three-dimensional lattice polytope, the second dilation is IDP.

\begin{proposition} \label{prop:IDP}
  Assume $(H,H)$ is IDP. If $D=kH$, $E=(k-1)H$ for $k \geq 2$, then
  $(D;E)$ has connected sections.
\end{proposition}

\begin{proof}
  As mentioned above, $(X,H)$ is projectively normal. Hence,
  we can decompose an arbitrary $\chi^m \in
  H^0(X,\OO(E))$ as $m=m_1+\ldots+m_{k-1}$. Furthermore, $m$ is
  connected in the section graph to $m'=m_1+\ldots+m_{k-2}+m_{k-1}'$ for arbitrary
  $\chi^{m'_{k-1}} \in H^0(X,\OO(H))$. By iterating over all indices,
  one obtains the connectedness of the section graph.
\end{proof}

\begin{corollary}
  Suppose $H$ is ample, $D=2kH$, $E=(2k-2)H$ for $k \ge 2$, then $(D;E)$
  has connected sections.
\end{corollary}

Let $\nefcone(X)$ denote the cone of nef divisors on $X$. 
\begin{proposition}\label{prop:cs}
  There is an ample $D_0$ on $X$ so that for every divisor $D \not\sim
  D_0$ with $[D] \in[ D_0] + \nefcone(X)$, the configuration
  $(D;D-D_0)$ has connected sections.
\end{proposition}
 To prove this proposition, we first prove two lemmas.
\begin{lemma} \label{lemma:ample-enough-GB}
	Fix a sequence \eqref{eq:modify} with corresponding matrix $A$. 
	There is an ample divisor $E_1$ such that for every divisor $E'$ satisfying $[E'] \in [E_1] + \nefcone(X)$, the set
  $(P(b') \cap \ZZ^p) - (P(b') \cap \ZZ^p)$
  is a Markov basis for $I_A$. Here, $P(b')$ is the polytope corresponding to $E'$.
\end{lemma}
\begin{proof}
Let $\generators$ be a finite Markov basis for $I_A$.
Identifying $\ker A$ with $M$, we may view $\generators$ as a finite subset of $M$. For any ample divisor $H$, there is an integer $k$ such that a translate of $P(kH)$ contains $\generators$. Taking $E_1=kH$, the claim follows from the fact that
\[P(E')=P(E_1)+P(E'-E_1).\] 
\end{proof}

\begin{lemma} \label{lemma:ample-enough-IDP}
	There is an ample divisor $E_2$ so that for all divisors $E',E$ with $[E'] \in [E_2] + \nefcone(X)$ 
	and $[E] \in \Nef(X)$ the pair $(E',E)$ is IDP.
\end{lemma}

\begin{proof}
The affine semigroup $\nefcone(X) \cap \Pic(X)$ has a finite generating
set $[D_1], \ldots, [D_r]$ by Gordon's lemma, where the $D_i$ are $T$-invariant  Cartier divisors. That is, every nef divisor
class can be represented by a non-negative integer combination of the
$[D_i]$.

Let $\delta := \max (C.D_i)$ be the highest degree of a
$\torus$-invariant curve $C$. Combinatorially, this is the longest edge length in the polytopes $P(D_i)$.
Choose $E_2 \in \Pic(X)$ ample enough to ensure $C.E_2 \ge 4\delta$
for all $C$, e.g., $E_2 = 4(D_1+\ldots+D_r)$.
Then for $[E'] \in [E_2] + \nefcone(X)$ and $[E] \in \nefcone(X)$ the pair
$(E',E)$ is IDP by \cite[Theorem 15 and Corollary 16]{convex-normal}.
\end{proof}

\begin{proof}[Proof of Proposition \ref{prop:cs}]
  We take $D_0$ to be $E_1+E_2$, where $E_1$ and $E_2$ are as in
  Lemmas \ref{lemma:ample-enough-GB} and
  \ref{lemma:ample-enough-IDP}.
  Then the pair $(D_0,D-D_0)$ is IDP by
  \ref{lemma:ample-enough-IDP}. Likewise, Lemma
  \ref{lemma:ample-enough-GB} together with Proposition
  \ref{prop:Markov} imply that the configuration $(D;D-D_0)$ has
  connected sections. The claim of the proposition follows.
\end{proof}

We conclude this section by considering an important class of toric
varieties. Recall that the root system of type $\cA_n$ is
\[
  \cA_n =  \{ \pm e_i \mid i = 1, \ldots, n \} \cup \{
  e_i-e_j \mid i,j = 1, \ldots, n,\, i \neq j \}\,\subset \ZZ^n.
\]
where $e_1,\ldots,e_n$ is the standard basis of $\ZZ^n$.
\begin{proposition}\label{prop:alcoved}
  Let $\fan$ be a fan with all rays generated by roots of type
  $\cA_n$, and let $X$ be the associated toric variety. Suppose $E$ and $E'$
  are nef divisors on $X$ with $E'$ big and set $D=E+E'$. Then $(D;E)$ has connected
  sections. 
\end{proposition}

\begin{proof}
The corresponding polytopes $P(E)$ and $P(E')$ are known as
type-$\cA$ polytopes or as alcoved polytopes. In this setting,
any pair $(E,E')$ is IDP \cite[Lemma~4.15]{triangulations}.
Further, every type-$\mathcal{A}$ polytope has a
canonical unimodular triangulation. The vertices of every simplex can be
ordered so that consecutive vertices differ by some $e^*_i$, an element of the basis dual to $\{e_i\}$
\cite[Theorem~3.3 and \S4.5]{triangulations}. This implies that
$\generators := \{ \pm e^*_i \mid i = 1, \ldots, n \}$ lifts to
generators of the toric ideal of $\mathcal{A}_n$. By
Proposition~\ref{prop:Markov} we see that $D,E$
has connected sections as soon as $D-E$ is big and nef.
\end{proof}

\section{Proofs of Main Results}\label{sec:proofs}
In this section, we will combine our lower bound on the genus (Theorem \ref{thm:z2}) with the discussion of \S \ref{sec:comb} to prove our main results Theorem \ref{thm:1} and Theorem \ref{thm:2}.
\begin{proof}[Proof of Theorem \ref{thm:1}]
	Since $D$ gives a projectively normal embedding, the pair $(D,D)$ is IDP. Proposition \ref{prop:IDP} then implies that for $m\geq 2$, the configuration $(mD;(m-1)D)$ has connected sections. The claim of the theorem then follows directly from Theorem \ref{thm:z2}.
\end{proof}

\begin{proof}[Proof of Theorem \ref{thm:2}]
  We first apply Proposition \ref{prop:cs} to obtain an ample divisor
  $D_0$ such that for any divisor $D \not\sim D_0$ whose class lies in
  $[D_0]+\Nef(X)$, $(D;D-D_0)$ has connected sections.
  Fix an ample class $H$.
  Let $H_0$ be any ample class such that $H_0-[D_0]$ and 
  $H_0-H-[D_0]+[K_X]$
   are ample, and all facets of the polytope
  corresponding to $H_0$ have at least two interior lattice points.
  For any divisor $D$ whose class lies in $H_0+\Nef(X)$, we thus still
  have that $(D;D-D_0)$ has connected sections.

  Let $S$ be any very general surface in $|D|$.
  Applying Theorem \ref{thm:z2} for the configuration $(D;D-D_0)$, we obtain
  for any integral curve $C$ not contained in the toric boundary of $S$,
  \[
    2g(C)-2\geq C.(D-D_0+K_X)\geq C.H.
  \]
  Thus, the only obstruction to the algebraic hyperbolicity of $S$ are the curves $C$ contained in the toric boundary. 

  For these curves, we may apply Lemma \ref{lemma:boundary} to see that they all have genus at least two. Since there are only finitely many of them, say $C_1,\ldots,C_k$, we may thus take 
  \[
    \epsilon=\min_{i=1,\ldots,k}\frac{1}{C_i.H}.
  \]
\end{proof}

\section{Examples}\label{sec:ex}
\begin{figure}[htb]
\centering
  \begin{tikzpicture}[scale=.5]
\draw [white,fill=lightgray] (5,3) -- (7,3) -- (7,7) -- (5,7) -- (5,4);
\draw [white,pattern color=black,pattern=crosshatch] (0,0) -- (7,0) -- (7,1) -- (3,1) -- (3,7) -- (0,7) -- (0,0);
		\foreach \x in {0,1,2,3,4,5,6,7} \foreach \y in {0,1,2,3,4,5,6,7} \draw [fill] (\x,\y) circle [radius=.08];
		\draw [thick,->] (0,0) -- (7.2,0);
		\draw [thick,->] (0,0) -- (0,7.2);
		\draw [thick,->] (5,3) -- (7.2,3);
		\draw [thick,->] (5,3) -- (5,7.2);
\node [below] at (7,0) {$a$};
\node [left] at (0,7) {$b$};
\node at (7,8) {hyperbolic};
\node at (1,8) {not hyperbolic};
	\end{tikzpicture}
	
	\caption{Algebraic hyperbolicity for very general surfaces of type $\CO(a,b)$ in $\PP^2\times\PP^1$}\label{fig:p2p1}
\end{figure}
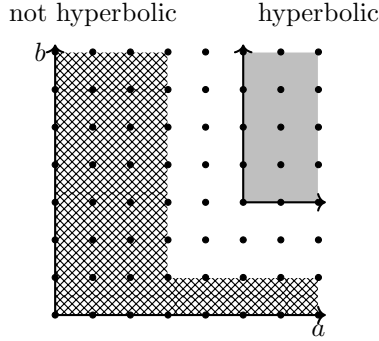

We now apply Theorem \ref{thm:z2} to obtain lower bounds on the genus
in some specific examples.
\begin{example}[$\PP^2\times \PP^1$]\label{ex:p2p1}
	Every line bundle on $X=\PP^2\times\PP^1$ is of the form 
\[\CO(a,b)=\pi_1^*(\CO_{\PP^2}(a))\otimes\pi_2^*(\CO_{\PP^1}(b)),\qquad a,b\in\ZZ
\]
where $\pi_1,\pi_2$ denote the projections of $\PP^2\times\PP^1$ onto
the first and second factors. By abuse of notation, we will write
$D=\CO(a,b)$ to mean that $D$ is a divisor whose associated line
bundle is isomorphic to $\CO(a,b)$.
Such a divisor $D$ is basepoint free if and only if $a,b\geq 0$ and
ample if and only if $a,b>0$. We also note that $K_X=\CO(-3,-2)$.
After fixing coordinates, we can assume that the fan $\Sigma$ associated to $\PP^2\times\PP^1$ has rays generated by $e_1,-e_2,e_2-e_1,e_3,-e_3$, so we are in a situation to apply Proposition \ref{prop:alcoved}.

Given $D=\CO(a,b)$ with $a,b\geq 2$, we can set $E=\CO(a-1,b-1)$. 
By Proposition \ref{prop:alcoved}, the configuration $(D;E)$ will have connected
sections. We thus obtain that for any curve $C$ not contained in the
toric boundary on a very general surface $S$ in $|D|$,
\[
	g\geq \frac{C.\CO(a-4,b-3)}{2}+1.
\]
In the case $a \ge 5$, $b \ge 4$ we choose $\epsilon=1$ for $H =
\CO(1,1)$ to get $2g-2 \ge \epsilon C.H$.

We use Lemma \ref{lemma:boundary} to analyze the curves contained in the boundary of
$S$, and notice that their geometric genera are exactly $(a-1)(b-1)$
and $(a-1)(a-2)/2$.
In particular, if $a\geq 4,b\geq 3$, $S$ contains no rational curves. 

We now assume that $a\geq 4,b\geq 3$.
By \cite[Theorem 4.2]{nl1}, it is straightforward to check that the
Noether-Lefschetz theorem holds for very general $S\in |D|$.
Thus, after tensoring with $\QQ$, any curve $C$ is rationally equivalent to the complete intersection of $S$ with a $\QQ$-divisor
of type $\CO(c,d)$, $c,d\geq 0$.
If $C$ is contained in the boundary, then we must have $(c,d)=(1,0)$ or $(0,1)$.
If $C$ is not contained in the boundary, an intersection number calculation yields 
\begin{align*}
	2g(C)-2\geq c(a(b-3)+(a-4)b)+da(a-4).
\end{align*}

On the other hand, the degree of such a curve $C$ with respect to the
polarization $H=\CO(1,1)$ is
\[
\deg C=C.H=c(a+b)+da.
\]

We claim that as long as $a\geq 5,b\geq 3$, $S$ is algebraically
hyperbolic. Indeed, we can take the constant $\epsilon$ to be
\[
  \epsilon=\frac{1}{a}.
\]
For $C$ not in the boundary we obtain
\[2g-2\geq cb+da\geq \frac{c(a+b)}{a}+d=\epsilon\cdot\deg C\]
as required.
For $C$ in the boundary, 
we have
\begin{align*}
  2g-2=2(a-1)(b-1)-2\geq \frac{a+b}{a}=\epsilon\cdot \deg C;\qquad(c,d)=(1,0);\\
  2g-2=(a-1)(a-2)-2\geq 1=\epsilon\cdot \deg C;\qquad(c,d)=(0,1)
\end{align*}
as required.

On the other hand, if $b\leq 1$ or $a\leq 3$ then $S$ contains curves
of genus zero or one, and cannot be algebraically hyperbolic. See
Figure \ref{fig:p2p1} for an illustration.
The only cases that remain open are when $a=4$ and $b\geq 2$, or $a\geq 4$ and $b=2$.
This proves Theorem \ref{thm:p2p1}.
\end{example}

\begin{example}[Example \ref{ex:p2p1} continued]\label{ex:cont}
	We continue our analysis of Example \ref{ex:p2p1} with a different configuration of divisors. As above, let $D=\CO(a,b)$ with $a,b\geq 1$. We set $E_1=\CO(a-1,b)$ and $E_2=\CO(a,b-1)$. It is straightforward to verify that $(D;E_1,E_2)$ has connected sections.

	Assuming that $a\geq 4,b\geq 3$, we may use the Noether-Lefschetz theorem  as above to assume that any curve $C\subset S$ is the intersection of $S$ with a divisor of type $\CO(c,d)$, $c,d\geq 0$. Assume that $(c,d)$ is not $(1,0)$ or $(0,1)$. Then $C$ is not contained in the boundary of $S$, and Theorem \ref{thm:z2} shows that
	\begin{align*}
g\geq \min \Big\{
	\frac{c(a(b-2)+b(a-4))+d(a(a-4))}{2}+1,\\
	\frac{c(a(b-3)+b(a-3))+d(a(a-3))}{2}+1
\Big\}.
	\end{align*}
As long as $c,d\geq 1$, this is better than the bound	
	\begin{align*}
g\geq 
	\frac{c(a(b-3)+b(a-4))+d(a(a-4))}{2}+1
	\end{align*}
	obtained in Example \ref{ex:p2p1} by just taking $E=\CO(a-1,b-1)$. Thus, if we know a bit more about $C$ than just the degree, we may obtain more refined lower bounds on the genus by taking configurations $(D;E_1,\ldots,E_\ell)$ involving multiple divisors.
\end{example}

\begin{example}[$\PP^1\times\PP^1\times\PP^1$]\label{ex:p1p1p1}
Every line bundle on $X=\PP^1\times\PP^1\times\PP^1$ is of the form 
\[\CO(a,b,c)=\pi_1^*(\CO_{\PP^1}(a))\otimes\pi_2^*(\CO_{\PP^1}(b))\otimes\pi_3^*(\CO_{\PP^1}(c)),\qquad a,b,c\in\ZZ
\]
where $\pi_1,\pi_2,\pi_3$ denote the projections onto the first, second, and third factors. By abuse of notation, we will write $D=\CO(a,b,c)$ to mean that $D$ is a divisor whose associated line bundle is isomorphic to $\CO(a,b,c)$.
Such a divisor $D$ is basepoint free if and only if $a,b,c\geq 0$ and ample if and only if $a,b,c>0$. We also note that $K_X=\CO(-2,-2,-2)$.
After fixing coordinates, we can assume that the fan $\Sigma$ associated to $\PP^1\times\PP^1\times\PP^1$ has rays generated by $e_1,-e_1,e_2,-e_2,e_3,-e_3$, so we are in a situation to apply Proposition \ref{prop:alcoved}.

Given $D=\CO(a,b,c)$ with $a,b,c\geq 2$, we can set 
\[E=\CO(a-1,b-1,c-1).\] 
The configuration $(D;E)$ has connected sections by Proposition \ref{prop:alcoved}.
We thus obtain that for any curve $C$ not contained in the toric boundary on a very general surface $S$ in $|D|$,
\[
	g\geq \frac{C.\CO(a-3,b-3,c-3)}{2}+1.
\]
We again use Lemma \ref{lemma:boundary} to analyze the curves contained in the boundary of $S$, and notice that their geometric genera are exactly $(a-1)(b-1)$, $(b-1)(c-1)$, and $(a-1)(c-1)$.
In particular, if $a,b,c\geq 3$, $S$ contains no rational curves. 

We now assume that $a\geq b\geq c\geq 3$.
By \cite[Theorem 4.2]{nl1}, it is straightforward to check that the Noether-Lefschetz theorem holds for very general $S\in |D|$.
Thus, after tensoring with $\QQ$, any curve $C$ is rationally equivalent to the complete intersection of $S$ with a $\QQ$-divisor
of type $\CO(d,e,f)$, $d,e,f\geq 0$.
If $C$ is contained in the boundary, then we must have $(d,e,f)=(1,0,0)$ or $(0,1,0)$ or $(0,0,1)$.
If $C$ is not contained in the boundary, an intersection number calculation yields
\begin{align*}
	2g(C)-2\geq d(b(c-3)+(b-3)c)+e(a(c-3)+(a-3)c)\\+f(a(b-3)+(a-3)b).
\end{align*}
On the other hand, the degree of such a curve $C$ with respect to the polarization $H=\CO(1,1,1)$ is
\[
\deg C=C.H=d(b+c)+e(a+c)+f(a+b).
\]

We claim that as long as $a\geq b > c\geq 3$, $S$ is algebraically hyperbolic. Indeed, we can take the constant $\epsilon$ to be \[
	\epsilon=\frac{1}{a+b}
\]
and obtain
\[2g-2\geq dc+ec+f(a+b)\geq \frac{d(b+c)}{a+b}+\frac{e(a+c)}{a+b}+f=\epsilon\cdot\deg C\]
for curves $C$ not in the boundary, as required.
For a curve $C$ in the boundary, we have
\begin{align*}
	2g-2=2(a-1)(b-1)-2\geq 1=\epsilon\cdot \deg C;\qquad(c,d)=(0,0,1);\\
	2g-2=2(a-1)(c-1)-2\geq \frac{a+c}{a+b}=\epsilon\cdot \deg C;\qquad(c,d)=(0,1,0);\\
	2g-2=2(b-1)(c-1)-2\geq \frac{b+c}{a+b}=\epsilon\cdot \deg C;\qquad(c,d)=(1,0,0)\\
\end{align*}
as required.

On the other hand, if $a$, $b$, or $c$ is less than two, then $S$ contains a rational curves and cannot be algebraically hyperbolic. Likewise, if $b=c=2$, $S$ contains an elliptic curve and similarly cannot be algebraically hyperbolic.
Assuming that $a\geq b \geq c$, the only cases that remain open are when $b>c=2$ or when $b=c=3$.
This proves Theorem \ref{thm:p1p1p1}.
\end{example}

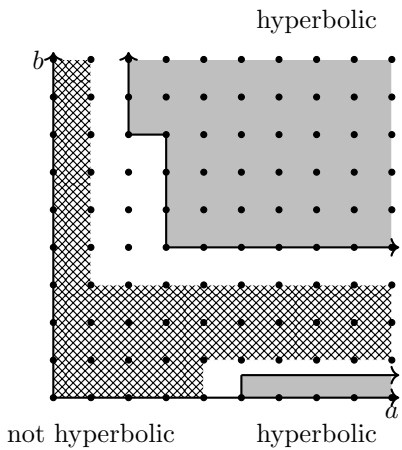
\begin{figure}[htb]
  \centering
	\begin{tikzpicture}[scale=.5]
\draw [white,fill=lightgray] (3,4) -- (9,4) -- (9,9) -- (2,9) -- (2,7) -- (3,7) -- (3,4);
\draw [white,fill=lightgray] (5,0) -- (9,0) -- (9,.6) -- (5,.6);
\draw [white,pattern color=black,pattern=crosshatch] (0,0) -- (4,0) -- (4,1) --(9,1) -- (9,3) -- (1,3) -- (1,9) -- (0,9) -- (0,0);
		\foreach \x in {0,1,2,3,4,5,6,7,8,9} \foreach \y in {0,1,2,3,4,5,6,7,8,9} \draw [fill] (\x,\y) circle [radius=.08];
		\draw [thick,->] (0,0) -- (9.2,0);
		\draw [thick,->] (0,0) -- (0,9.2);
		\draw [thick,->] (5,.6) -- (9.2,.6);
		\draw [thick] (5,0) -- (5,.6);
		\draw [thick] (2,7) -- (3,7) --(3,4);
		\draw [thick,->] (2,7) -- (2,9.2);
		\draw [thick,->] (3,4) -- (9.2,4);
\node [below] at (9,0) {$a$};
\node [left] at (0,9) {$b$};
\node at (7,10) {hyperbolic};
\node at (1,-1) {not hyperbolic};
\node at (7,-1) {hyperbolic};
	\end{tikzpicture}

	\caption{Algebraic hyperbolicity for very general surfaces of type $aH+b(H-E)$ in the blowup of $\PP^3$}\label{fig:blp32}
\end{figure}

\begin{example}[Blowup of $\PP^3$ at a point]\label{ex:bl}
	Let $X$ be the blowup of $\PP^3$ at a point. We take $H$ to be the pullback of a hyperplane, and $E$ the exceptional divisor. These two divisors generate the Picard group of $X$. The nef cone is generated by $H$ and $H-E$, whereas the effective cone is generated by $E$ and $H-E$. A canonical divisor on $X$ is $K_X=-2H-2(H-E)$.
We also recall the following intersection products:
\begin{align*}
E.E.E=H.H.H=1\\
E.E.H=E.H.H=0
\end{align*}
After fixing coordinates, we can assume that the fan $\Sigma$ associated to $X$ has rays generated by $e_1,e_2,e_2-e_3,e_3-e_1,e_1-e_3$, so we are in a situation to apply Proposition \ref{prop:alcoved}.

\begin{figure}[htb]
  \centering
	\begin{tikzpicture}
\draw (5,0) -- (0,0) -- (0,5);
\draw [dashed] (5,0) -- (0,5);
\draw (0,0) -- (1,2);
\draw (1,2) -- (3.5,2) -- (1,4.5) -- (1,2);
\draw (3.5,2) -- (5,0);
\draw (1,4.5) -- (0,5);
\node [below] at (2.5,0) {$a+b$};
\node [left] at (0,2.5) {$a+b$};
\node [below right] at (.5,1) {$a$};
\node [below] at (2,2) {$b$};
\node [left] at (1,3) {$b$};
	\end{tikzpicture}

	\caption{The polytope $P(D)$ for the blowup of $\PP^3$}\label{fig:blp3}
\end{figure}
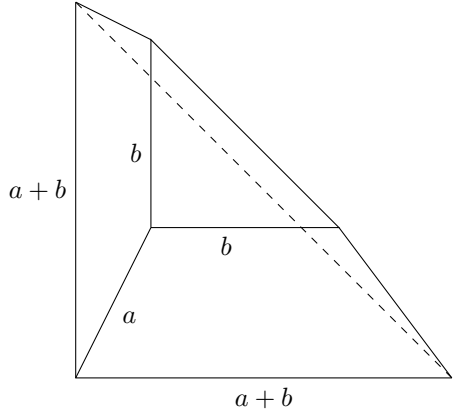

	Consider $D=aH+b(H-E)$ for $a\geq 1,b\geq 0$. The polytope corresponding to $D$ is pictured in Figure \ref{fig:blp3} (although we use different coordinates than chosen above).
	Set $E_1=(a-1)H+b(H-E)$.
	By Proposition \ref{prop:alcoved}, the configuration $(D;E_1)$ has connected sections.
	We obtain that for any curve $C$ not contained in the toric boundary on a very general surface $S$ in $|D|$,
\[
	g\geq \frac{C.(a-3)H+(b-2)(H-E))}{2}+1.
\]
The curves contained in the boundary of $S$ have geometric genera 
\[
	\frac{(a+b-1)(a+b-2)}{2},\frac{(b-1)(b-2)}{2},\frac{(a+b-1)(a+b-2)-b(b-1)}{2}.
\]

We now assume that $a,b\geq 2$.
The variety $X$ is smooth of Picard rank two, so the hypotheses of \cite[Theorem 4.2]{nl1} are fulfilled for $D$, see e.g.~\cite[Section 3.2]{nl2}.
Thus,
after tensoring with $\QQ$, any curve $C$ is rationally equivalent to the complete intersection of $S$ with a $\QQ$-divisor
 of type $c(H-E)+dE$, $c,d\geq 0$.
If $C$ is contained in the boundary, then we must have $(c,d)=(1,1)$ or $(0,1)$ or $(1,0)$.
If $C$ is not contained in the boundary, an intersection number calculation yields
\begin{align*}
	2g(C)-2\geq c(a(a-3)+a(b-2)+b(a-3))+d(b(b-2)).
\end{align*}
On the other hand, the degree of such a curve $C$ with respect to the polarization $H+(H-E)$ is
\[
\deg C=c(2a+b)+db.
\]

We claim that as long as $a\geq 3$ and $b\geq 4$, $S$ is algebraically hyperbolic. 
Likewise, if $a=2$ and $b\geq 7$, $S$ is algebraically hyperbolic.
Indeed, in both cases we can take the constant $\epsilon$ to be \[
	\epsilon=\frac{1}{4a+2b};
\]
details are left to the reader.

On the other hand, as long as $b\geq 1$, if or $ a<2$ or $b<4$ then $S$ contains curves of genus zero or one in its boundary, and cannot be algebraically hyperbolic. 
If $b=0$, then $S$ is just a very general surface of degree $a$ in $\PP^3$. This is algebraically hyperbolic by \cite{coskun} if and only if $a\geq 5$; our methods suffice to show hyperbolicity as long as $a\geq 6$. 
See Figure \ref{fig:blp32} for an illustration of when $S$ is algebraically hyperbolic.
The only cases that remain open are when $a=2$ and $b=4,5,6$;
This proves Theorem \ref{thm:blp3}.
\end{example}

We can also use our techniques in the case of non-Gorenstein singularities, although a bit more care is required:
\begin{example}[Weighted projective space $\PP(1,1,1,n)$]\label{ex:wps}
	We consider the weighted projective space $X=\PP(1,1,1,n)$. This has an isolated, non-Gorenstein singularity. Let $H$ be an ample (Cartier) generator of the Picard group of $X$; the corresponding sheaf is often denoted $\CO(n)$. We will consider curves $C$ on a very general surface $S\in |D|$, where $D=mH$ for some $m\geq 1$.

	To apply our results, we need to resolve singularities. Let $\pi:\widetilde{X}\to X$ be the blowup of $X$ at the singular point. The result is a smooth toric variety $\widetilde{X}$ of Picard number two; the Picard group (and nef cone) are generated by the pullback of $H$ (which we also denote by $H$) and a divisor $F$ satisfying $nF=H-E$ in the Picard group. Here, $E$ is the exceptional divisor of the blowup. A canonical divisor on $\widetilde{X}$ is given by $K_{\widetilde{X}}=-2H-(3-n)F$. For the special case $n=2$, we make the important observation that for any curve $\widetilde{C}$ in $\widetilde{X}$, if $\widetilde{C}.(H-2F)<0$, then $\widetilde{C}$ must be contained in $E$, and hence is contracted by $\pi$. Indeed, in this case, $H-2F$ is equivalent to $E$.

	Let $S$ be a very general surface in $|mH|$ on $X$, and $\widetilde{S}$ the pullback to $\widetilde{X}$; this is a very general surface in $|mH|$ on $\widetilde{X}$. Since $H$ gives a projectively normal embedding of $X$, we may apply Proposition \ref{prop:IDP} to conclude that the pair $D=mH,E_1=(m-1)H$ has connected sections. Let $C$ be any curve in $S$ not contained in the toric boundary, and $\widetilde{C}$ its preimage in $\widetilde S$. Then 
	\[
		g(C)=g(\widetilde{C})\geq \frac{\widetilde{C}.((m-3)H+(n-3)F)}{2}+1.
	\]

For $n\geq 3$, we obtain
	\[
		g(C)=g(\widetilde{C})\geq \frac{{C}.(m-3)H}{2}+1.
	\]
by the projection formula.
Suppose instead that $n=2$.
	By construction, $\widetilde{C}$ is not contracted by $\pi$, so $\widetilde{C}.(H-2F)\geq 0$. 
We thus have
\begin{align*}
	g(C)\geq \frac{\widetilde{C}.((m-4)H+(H-2F)+F)}{2}+1\geq \frac{\widetilde{C}.(m-4)H}{2}+1\\
=\frac{C.(m-4)H}{2}+1
\end{align*}
by the projection formula.

On the other hand, for any $n\geq 2$, a curve in the toric boundary of $S$ will have genus $(nm-1)(nm-2)/2$ or $(m-1)(nm-2)/2$. It follows that $X$ is algebraically hyperbolic if $n\geq 3$ and $m\geq 4$, or $n=2$ and $m\geq 5$. If $m=1$ or $m=n=2$, $X$ has curves of genus less than two, so is not algebraically hyperbolic. The open cases are $n=2,m=3,4$ and $n\geq3,m=2,3$.
This proves Theorem \ref{thm:wps}.
\end{example}

These examples lead to the following question:
\begin{question}\label{q:hyper}
Are very general surfaces of the following types algebraically hyperbolic?
\begin{enumerate}
	\item $\CO(a,b)$ in $\PP^2\times\PP^1$ with $a=4,b\geq 3$ or $a\geq 4,b=3$;
	\item $\CO(a,b,c)$ in $\PP^1\times\PP^1\times\PP^1$ with $a\geq b = c =3$ or $b>c=2$;
	\item $2H+b(H-E)$ in the blowup of $\PP^3$ at a point for $b=4,5,6$;
	\item $\CO(6)$ and $\CO(8)$ in weighted projective space $\PP(1,1,1,2)$;
	\item $\CO(2n)$ and $\CO(3n)$ in weighted projective space $\PP(1,1,1,n)$ for $n\geq 3$.
\end{enumerate}
\end{question}
\noindent This question will be resolved in forthcoming work by Coskun and Riedl \cite{forth}.

\end{document}